\documentclass[10pt,leqno]{amsart}
\usepackage{graphicx}
\usepackage{indentfirst,csquotes}

\usepackage{amssymb,amsthm,amsmath}
\usepackage{xcolor,paralist,hyperref,fancyhdr,etoolbox}
\newtheorem{th.}{Theorem}[]
\newtheorem{de.}[th.]{Definition}
\newtheorem{ex.}[th.]{Example}
\newtheorem{le.}[th.]{Lemma}
\newtheorem{re.}[th.]{Remark}
\newtheorem{pr.}[th.]{Proposition}
\newtheorem{co.}[th.]{Corollary}
\newtheorem{conj.}[th.]{Conjecture}

\hypersetup{ colorlinks=true, linkcolor=black, filecolor=black, urlcolor=black }

\usepackage{lipsum}

\begin{document}
\title{An integration of Lie-Leibniz triples} 
\author{Ryo Hayami}
\date{\today}
\address{Nagano University}
\email{ryo-hayami@nagano.ac.jp}

\let\thefootnote\relax
\footnotemark 

\begin{abstract}
In this paper, we introduce the group version of a Lie-Leibniz triple, which we call a Lie group-rack triple. We define a Lie group-rack triple whose tangent structure is a Lie-Leibniz triple, which is a generalization of an augmented Lie rack whose tangent structure is an augmented Leibniz algebra. We show that any finite-dimensional Lie-Leibniz triple can be integrated to a local Lie group-rack triple by generalizing the integration procedure of an augmented Leibniz algebra into an augmented Lie rack.   
\end{abstract}
\maketitle

\tableofcontents

\bigskip

\section{introduction}

Lie-Leibniz triples appear when physicists formulate gauged models of supergravity\cite{deWit:2005hv,deWit:2008gc}. In gauged supergravity, for gauging of a subgroup of a global symmetry group, we have to add a hierarchy of higher-form fields to maintain the covariance of the theory under the global group. This hierarchy is called a tensor hierarchy. This hierarchy is based on a Lie algebra $\mathfrak{g}$, a choice of representation $V$, and a linear map $\theta:V\rightarrow\mathfrak{g}$ satisfying some compatibility conditions, which is called an embedding tensor. Owing to the conditions of $\theta$, $V$ has a structure of a Leibniz algebra, a vector space with a bilinear bracket $[,]$ which satisfies 
\begin{equation}
[u,[v,w]]=[[u,v],w]+[v,[u,w]]
\end{equation}
for all $u,v,w\in V$. The mathematical formulation for the above $(\mathfrak{g},V,\theta)$ is a Lie-Leibniz triple.  

Let $\mathfrak{g},[,]_{\mathfrak{g}}$ be a Lie algebra, and $V,[,]$ be a Leibniz algebra. Assume that there is a Lie algebra action of $\mathfrak{g}$ on $V$. An embedding tensor is a linear map $\theta:V\rightarrow \mathfrak{g}$ which satisfies two compatibility conditions. The triple $(\mathfrak{g},V,\theta)$ is called a Lie-Leibniz triple. This triple gives rise to a differential-graded Lie algebra (dgLa for short). This dgLa is the mathematical formulation of the tensor hierarchy\cite{Palmkvist:2013vya,Lavau:2017tvi}.

The aim of this paper is to give a group version of Lie-Leibniz triples and to develop some integration theory of these triples. In other words, we try to find some kind of a differential-geometric structure whose tangent structure is a given Lie-Leibniz triple.  

It is well known that the tangent space of a Lie group at the unit element has the structure of a Lie algebra. Lie's third theorem states that for any finite-dimensional real Lie algebra $\mathfrak{g}$ there is a Lie group $G$ such that the tangent Lie algebra of $G$ is isomorphic to $\mathfrak{g}$. The question of finding a generalization of Lie's third theorem for Leibniz algebras was formulated in \cite{loday1993version}. In \cite{kinyon2007leibniz}, a Lie rack was introduced as an answer to the question.

The definition of a Lie rack is as follows. 

\begin{de.}

A set $X$ with a binary operation $\triangleright:X\times X \rightarrow X$ is a rack if it satisfies 
\begin{equation}
x\triangleright(y\triangleright z)=(x\triangleright y)\triangleright(x\triangleright z)
\end{equation}
for all $x,y,z\in X$ (this relation is called self-distributivity) and the map $x\triangleright-:X\rightarrow X$ is bijective. A rack $(X,\triangleright)$ with a distinguished point $e\in X$ is a pointed rack if $e\triangleright x=x, x\triangleright e=e$ for all $x\in X$. When $X$ is a smooth manifold and $\triangleright$ is smooth, $(X,\triangleright,e)$ is a Lie rack.

\end{de.}

For a (Lie) group $G$, $(G,\triangleright,e_{G})$ is a (Lie) rack where $g_{1}\triangleright g_{2}=g_{1}g_{2}g_{1}^{-1}$ and $e_{G}$ is the unit element of $G$. This rack is called a conjugation (Lie) rack of $G$.

The tangent space of a Lie rack $X$ at $e$ with the infinitesimal bracket of the operation $\triangleright$ satisfies the conditions of a Leibniz algebra. Therefore, Lie racks can be seen as integrated objects of Leibniz algebras.

In \cite{kinyon2007leibniz}, the way to integrate Leibniz algebras into Lie racks was proposed, but this is not a generalization of the way to integrate Lie algebras into the conjugation rack of the associated Lie groups. In \cite{covez2010integration}, the local integration of Leibniz algebras into local Lie racks, which is a generalization of the integration of Lie algebras into local Lie groups, was introduced. The global integration of general Leibniz algebras into Lie racks which reduces in the case of Lie algebras to the ordinary integration into Lie groups was constructed in \cite{bordemann2017global}. The point of this integration is that a Leibniz algebra with an additional structure, which is called an augmented Leibniz algebra, is integrated into a Lie rack with an additional structure, which is called an augmented Lie rack.

In fact, an augmented Leibniz algebra is equivalent to a class of Lie-Leibniz triples, which is called a strict Lie-Leibniz triple. Therefore, in terms of Lie-Leibniz triples, the result in \cite{bordemann2017global} is that strict Lie-Leibniz triples can be integrated into augmented Lie racks. In this paper, we introduce a group version of Lie-Leibniz triples, which we call Lie group-rack triples, and show that a class of Lie group-rack triples which we call a relaxed augmented Lie rack can be seen as an integrated object of a (not necessarily strict) Lie-Leibniz triple. Moreover, we give an explicit local integration method into local Lie group-rack triples for finite-dimensional Lie-Leibniz triples. 

First, we recall the definition of Lie-Leibniz triples and study the relation to augmented Leibniz algebras. In particular, we show that Lie-Leibniz triples are equivalent to a little generalization of augmented Leibniz algebras, which we call relaxed augmented Leibniz algebras. Next, we define a Lie group-rack triple and study its properties. We also define a relaxed augmented Lie rack as an example of Lie group-rack triples and check that its tangent structure is equivalent to a Lie-Leibniz triple. Moreover, we give a local integration into a local relaxed augmented Lie rack for a given Lie-Leibniz triple. The integration procedure generalizes the local structure of the integration introduced in \cite{bordemann2017global}.

\section{Lie-Leibniz triples and augmented Leibniz algebras}

In this section, we recall the basics of the Lie-Leibniz triples and study the relation with the augmented Leibniz algebras.

The definition of a Lie-Leibniz triple is as follows.

\begin{de.}

A Lie-Leibniz triple is a triple $(\mathfrak{g},V,\theta)$ where $(\mathfrak{g},[,]_{\mathfrak{g}})$ is a Lie algebra, $(V,[,])$ is a Leibniz algebra which is also a $\mathfrak{g}$-module and $\theta:V\rightarrow\mathfrak{g}$ is a linear map called the embedding tensor.

These structures satisfy the following two conditions for all $u,v\in V$.

\begin{description}

\item[linear constraint]$[u,v]=\theta(u)\cdot v$. 

\item[quadratic constraint] $\theta[u,v]=[\theta(u),\theta(v)]_{\mathfrak{g}}$.

\end{description}

\end{de.}

\begin{de.}

A morphism between two Lie-Leibniz triples $(\mathfrak{g},V,\theta)$ and $(\mathfrak{g}',V',\theta')$ is a pair $(\varphi,\psi)$ where $\varphi:\mathfrak{g}\rightarrow\mathfrak{g}'$ is a Lie algebra morphism and $\psi:V\rightarrow V'$ is a morphism of vector spaces such that

\begin{equation}
\theta'(\psi(v))=\varphi(\theta(v))
\end{equation}
\begin{equation}
\psi(av)=\varphi(a)\cdot\psi(v)
\end{equation}
holds for every $a\in\mathfrak{g}$ and $x\in V$.

Two Lie-Leibniz algebras are isomorphic if there is a morphism $(\varphi,\psi)$ between the Lie-Leibniz algebras where $\varphi$ and $\psi$ are isomorphisms.
\end{de.}

From the definition, we can see that $\psi$ becomes a Leibniz algebra morphism. Indeed, for all $u,v\in V$,
\begin{align*}
[\psi(u),\psi(v)]&=\theta'(\psi(u))\cdot\psi(v)\\
&=\varphi(\theta(u))\cdot\psi(v)\\
&=\psi(\theta(u)\cdot v)=\psi([u,v]).\\
\end{align*}

Next, we recall the definition of augmented Leibniz algebras.

\begin{de.}

Let $\mathfrak{g}$ be a Lie algebra and act on a vector space $V$. We say that $V$ together with a linear map $\theta:V\rightarrow \mathfrak{g}$ is an augmented Leibniz algebra when it satisfies
\begin{equation}
\theta(g\cdot v)=[g,\theta(v)]_{\mathfrak{g}}
\end{equation}
for all $g\in \mathfrak{g}$ and all $v\in V$. 

\end{de.}

Let $V$ be an augmented Leibniz algebra with $\theta:V\rightarrow \mathfrak{g}$. Then $V$ has a Leibniz bracket given by
\begin{equation}
[u,v]:=\theta(u)\cdot v.
\end{equation}
We can see that $(V,[,])$ is a Leibniz algebra by the following calculation.
\begin{align*}
[u,[v,w]]&=\theta(u)\cdot(\theta(v)\cdot w)\\
&=\theta(u)\cdot(\theta(v)\cdot w)-\theta(v)\cdot(\theta(u)\cdot w)+\theta(v)\cdot(\theta(u)\cdot w)\\
&=[\theta(u),\theta(v)]_{\mathfrak{g}}\cdot w+\theta(v)\cdot(\theta(u)\cdot w)\\
&=\theta(\theta(u)\cdot v)\cdot w+\theta(v)\cdot(\theta(u)\cdot w)\\
&=[[u,v],w]+[v,[u,w]].
\end{align*}
Moreover, we can see that the $\mathfrak{g}$-action gives a derivation of $V$ as Leibniz algebras, and that $\mathrm{Im}(\theta)$ is an ideal of $\mathfrak{g}$, as follows. For all $g\in\mathfrak{g}$ and $u,v\in V$, 
\[
[gu,v]=\theta(gu)\cdot v=[g,\theta(u)]_{\mathfrak{g}}\cdot v=g\cdot(\theta(u)v)-\theta(u)\cdot(gv)=g[u,v]-[u,gv],
\]
\[
[g,\theta(v)]_{\mathfrak{g}}=\theta(g\cdot v)\in\mathrm{Im}(\theta).
\]

For an augmented Leibniz algebra, $(\mathfrak{g},V,\theta)$ becomes a Lie-Leibniz triple. In \cite{chen2025dg}, the equivalent structure is called a Loday-Pirashvili module.

We define a little generalization of an augmented Leibniz algebra by relaxing the condition of $\theta$ by replacing "for all $g\in \mathfrak{g}$" by "for all $g\in \mathfrak{g}'$" where $\mathfrak{g}'$ is a Lie subalgebra of $\mathfrak{g}$ including $\mathrm{Im}(\theta)$. 

\begin{de.}

Let $\mathfrak{g}$ be a Lie algebra and act on a vector space $V$. We say that $V$ together with a linear mapping $\theta:V\rightarrow \mathfrak{g}$ and a Lie subalgebra $\mathfrak{g}'\subset\mathfrak{g}$ including $\mathrm{Im}(\theta)$ is a relaxed augmented Leibniz algebra when it satisfies
\begin{equation}
\theta(g\cdot v)=[g,\theta(v)]_{\mathfrak{g}}
\end{equation}
for all $g\in \mathfrak{g}'$ and all $v\in V$. 

\end{de.}

Let $V$ be a relaxed augmented Leibniz algebra with $\theta:V\rightarrow \mathfrak{g}$ and $\mathfrak{g}'\subset\mathfrak{g}$. Then $V$ also has a Leibniz bracket given by
\begin{equation}
[u, v]:=[\theta(u), v]_{\mathfrak{g}}.
\end{equation}
and $(\mathfrak{g},V,\theta)$ also becomes a Lie-Leibniz triple. We can check the Leibniz identity in the same way for an ordinary augmented Leibniz algebra. Moreover, $(\mathfrak{g}',V,\theta)$ becomes an augmented Leibniz algebra, viewing $V$ as $\mathfrak{g}'$-module and $\theta$ as a map from $V$ to $\mathfrak{g}'$, and $\mathfrak{g}$-action of $h\in\mathfrak{h}$ gives a derivation of $V$ as Leibniz algebras and $\mathrm{Im}(\theta)$ is an ideal of $\mathfrak{h}$. However, the $\mathfrak{g}$-action of $g\in\mathfrak{g}$ does not have to give a derivation of $V$ as Leibniz algebras when $g\notin\mathfrak{h}$, and $\mathrm{Im}(\theta)$ is not necessarily an ideal of $\mathfrak{g}$.

Let $(\mathfrak{g},V,\theta)$ be a Lie-Leibniz triple. For $a\in\mathfrak{g}$, we define a map $a\theta\in\mathrm{Hom}(V,\mathfrak{g})$ by
\begin{equation}
a\theta:V\ni v\mapsto [a,\theta(v)]_{\mathfrak{g}}-\theta(a\cdot v)\in\mathfrak{g}.
\end{equation}

 A Lie-Leibniz triple is called strict if $a\theta=0$ for all $a\in\mathfrak{g}$. Any Lie-Leibniz triple from an augmented Leibniz algebra is strict. Conversely, a strict Lie-Leibniz triple $(\mathfrak{g},V,\theta)$ has the structure of an augmented Leibniz algebra.

In general, for a Lie-Leibniz triple $(\mathfrak{g},V,\theta)$, $a\theta=0$ is satisfied for all $a\in\mathrm{Im}(\theta)$ but not for all $a\in\mathfrak{g}$. Take a Lie subalgebra $\mathfrak{g'}$ of $\mathfrak{g}$ such that $a\theta=0$ is satisfied for all $a\in\mathfrak{g'}$ (for example, $\mathrm{Im}(\theta)$ is itself a Lie subalgebra of $\mathfrak{g}$). Then $(\mathfrak{g},V,\theta,\mathfrak{g'})$ becomes a relaxed augmented Leibniz algebra.

An important example of augmented Leibniz algebras is a Lie algebra crossed module. We refer to \cite{wagemann2021crossed} for the basics of crossed modules.

\begin{de.}

A Lie algebra crossed module is a 4-tuple $(\mathfrak{m},\mathfrak{n},\mu:\mathfrak{m}\rightarrow\mathfrak{n},\eta:\mathfrak{n}\rightarrow\mathrm{End}(\mathfrak{m}))$, where $\mathfrak{m}$ and $\mathfrak{n}$ are Lie algebras, $\mu$ is a homomorphism of vector spaces, and $\eta$ is a Lie algebra action of $\mathfrak{n}$ on $\mathfrak{m}$, which satisfies the following conditions for all $n\in\mathfrak{n}$ and $m,m'\in\mathfrak{m}$.
\begin{enumerate}

\item$\mu(\eta(n)(m))=[n,\mu(m)]_{\mathfrak{n}}$.

\item$\eta(\mu(m))(m')=[m,m']_{\mathfrak{m}}$.

\end{enumerate}

\end{de.}

We can check that $\mu$ becomes a homomorphism of Lie algebras and that $\eta(n)\in\mathrm{Der}(\mathfrak{m})$ as follows.
\[
\mu([m,m']_{\mathfrak{m}})=\mu(\eta(\mu(m))(m'))=[\mu(m),\mu(m')],
\]
\begin{align*}
[\eta(n)(m),m']_{\mathfrak{m}}&=\eta(\mu(\eta(n)(m)))(m')=\eta([n,\mu(m)]_{\mathfrak{n}})(m')\\
&=\eta(n)(\eta(\mu(m))(m'))-\eta(\mu(m))(\eta(n)(m'))=\eta(n)([m,m']_{\mathfrak{m}})-[m,\eta(n)(m')]_{\mathfrak{m}}.
\end{align*}

Viewing $\eta:\mathfrak{n}\rightarrow\mathfrak{der(m)}$ as an $\mathfrak{n}$-module structure on $\mathfrak{m}$, $(\mathfrak{n},\mathfrak{m},\mu)$ becomes an augmented Leibniz algebra.

We can relax the condition in a similar way that augmented Leibniz algebras are relaxed.

\begin{de.}

A relaxed Lie algebra crossed module is a 5-tuple $(\mathfrak{m},\mathfrak{n},\mathfrak{n}',\mu:\mathfrak{m}\rightarrow\mathfrak{n},\eta:\mathfrak{n}\rightarrow \mathrm{End}(\mathfrak{m}))$, where $\mathfrak{n}'\subset\mathfrak{n}$ is a Lie subalgebra of $\mathfrak{n}$ including $\mathrm{Im}(\mu)$, such that the two conditions of Lie algebra crossed module holds for all $n\in \mathfrak{n}',m,m'\in \mathfrak{m}$.

\end{de.}

We can check that $\mu$ becomes a homomorphism of Lie algebras completely analogously to the original case, because $\mathrm{Im}(\mu)\subset\mathfrak{n}'$, but $\eta(n)\in\mathrm{Der}(\mathfrak{m})$ does not need to be satisfied when $n\notin\mathfrak{n}'$. Viewing $\eta:\mathfrak{n}\rightarrow\mathrm{End}\mathfrak{(m)}$ as an $\mathfrak{n}$-module structure on $\mathfrak{m}$, $(\mathfrak{n},\mathfrak{m},\mu:\mathfrak{m}\rightarrow\mathfrak{n},\mathfrak{n}')$ is a relaxed augmented Leibniz algebra and hence $(\mathfrak{n},\mathfrak{m},\mu)$ is a Lie-Leibniz triple. Lie-Leibniz triples of this form appeared in Proposition 13 in \cite{Lavau:2020pwa}.



\section{Lie group-rack triples and augmented Lie racks}

In this section, we define Lie group-rack triples and study their properties. We assume that all manifolds and Lie groups in this section are finite-dimensional.

\begin{de.}

A group-rack triple is a triple $(G,X,\Theta)$ where $G$ is a group with a unit element $e_{G}$, $(X,\triangleright,e)$ is a pointed rack where $X$ is a $G$-set and $\Theta:X\rightarrow G$ is a map that maps $e\in X$ to a unit element $e_{G}$ and that satisfies the following two conditions for all $x,y\in X$. 

\begin{description}

\item[linear constraint]:$x\triangleright y=\Theta(x)\cdot y$. 

\item[quadratic constraint]: $\Theta(x\triangleright y)=\Theta(x)\Theta(y)\Theta^{-1}(x)$.

\end{description}

A Lie group-rack triple is a group-rack triple $(G,X,\Theta)$ where $G$ is a Lie group, $(X,\triangleright,e)$ is a Lie rack on which $G$ acts smoothly and $\Theta:X\rightarrow G$ is a smooth mapping. 

\end{de.}

\begin{re.}

The equivalent notion of a group-rack triple had already appeared in 2.24 of \cite{das2024averaging}. There this notion is called a relative averaging operator.

\end{re.}


A morphism of Lie group-rack triples is defined as follows.

\begin{de.}

A morphism between two Lie group-rack triples $(G,X,\Theta)$ and $(G',X',\Theta')$ is a pair $(\Phi,\Psi)$ where $\Phi:G\rightarrow G'$ is a morphism of Lie groups and $\Psi:X\rightarrow X'$ is a morphism of smooth manifolds such that
\begin{equation}
\Psi(e_{X})=e_{X'}    
\end{equation}
\begin{equation}
\Theta'(\Psi(x))=\Phi(\Theta(x))
\end{equation}
\begin{equation}
\Psi(g\cdot x)=\Phi(g)\cdot\Psi(x)
\end{equation}
for every $g\in G,x\in X$.
\end{de.}

We can check that $\Psi$ is a morphism of Lie racks as follows:

\begin{align*}
\Psi(x\triangleright y)&=\Psi(\Theta(x)\cdot y)\\
&=\Phi(\Theta(x))\cdot\Psi(y)\\
&=\Theta'(\Psi(x))\cdot\Psi(y)\\
&=\Psi(x)\triangleright\Psi(y).
\end{align*}

Let $(G,X,\Theta)$ be a group-rack triple. For $g\in G$, we define a map $g\Theta\in C^{\infty}(X,G)$ by
\begin{equation}
g\Theta:X\ni x\mapsto (g\Theta(x)g^{-1})\Theta(gx)^{-1}\in G.
\end{equation}
For any group-rack triple, $g\Theta=e_{G}$ holds for all $g\in\mathrm{Im}(\Theta)$. When $g\Theta=e_{G}$ for all $g\in G$, we say that the group-rack triple is strict.

An important fact about a strict group-rack triple is that this is equivalent to an augmented rack which was introduced in \cite{fenn1992racks}. The Lie version of an augmented rack is called an augmented Lie rack.

\begin{de.}

Let $G$ be a Lie group and act on a manifold $X$ with a point $e\in X$ such that $g\cdot e=e$ for all $g\in G$. We say that $X$ with a smooth map $\Theta:X\rightarrow G$ is an augmented Lie rack if
\begin{equation}
\Theta(e)=e_{G}    
\end{equation}
and
\begin{equation}\label{ar}
\Theta(g\cdot x)=g\Theta(x)g^{-1}
\end{equation}
holds for all $g\in G$ and all $x\in X$. 

\end{de.}

Let $X$ be an augmented Lie rack with $\Theta:X\rightarrow G$ and $e\in X$. Then $X$ has a rack product given by
\begin{equation}
x\triangleright y:=\Theta(x)\cdot y.
\end{equation}
for all $x,y\in X$. We can check that $(X,\triangleright)$ becomes a Lie rack by the following calculations.
\begin{align*}
x\triangleright(y\triangleright z)&=\Theta(x)\cdot(\Theta(y)\cdot z)\\
&=\Theta(x)\Theta(y)\Theta(x)^{-1}\Theta(x)\cdot z\\
&=(\Theta(\Theta(x)\cdot y))\cdot(\Theta(x)\cdot z)=(x\triangleright y)\triangleright(y\triangleright z),
\end{align*}

\begin{equation*}
e\triangleright x=e_{G}\cdot x=x,\ x\triangleright e=\Theta(x)\cdot e=e.
\end{equation*}

If $(G,X,p)$ is an augmented Lie rack, then it becomes a strict Lie group-rack triple and the converse also holds.

In \cite{bordemann2017global}, the authors presented a global integration procedure for any real finite-dimensional augmented Leibniz algebra into an augmented Lie rack. This means that any real finite-dimensional strict Lie-Leibniz triple can be integrated into a strict Lie group-rack triple.

We define a little generalization of an augmented Lie rack by relaxing the condition (\ref{ar}) by replacing "for all $g\in G$" by "for all $g\in G'$" where $G'$ is a Lie subgroup of $G$ containing $\mathrm{Im}(\Theta)$.

\begin{de.}
Let $G$ be a Lie group and act on a manifold $X$ with a distiguished point $e\in X$. We say that $X$ with a smooth map $\Theta:X\rightarrow G$ and a Lie subgroup $G'\subset G$ including $\mathrm{Im}(\Theta)$ is a relaxed augmented Lie rack if
\begin{equation}
\Theta(e)=e_{G}    
\end{equation}
and
\begin{equation}
\Theta(g\cdot x)=g\Theta(x)g^{-1}
\end{equation}
holds for all $g\in G'$ and all $x\in X$. 

\end{de.}

Let $(X,G,\Theta:X\rightarrow G,G')$ be a relaxed augmented Lie rack. Then $X$ has a rack product given by
\begin{equation}
x\triangleright y:=\Theta(x)\cdot y.
\end{equation}
The axioms of a Lie rack can be checked in the same way for the strict(non-relaxed) case. When $G'=G$, the relaxed augmented Lie rack reduces to an ordinary augmented Lie rack. For a relaxed augmented Lie rack, $(G,X,\Theta)$ becomes a Lie group-rack triple. Note that all the above discussions about relaxed augmented Lie racks can be generalized to the non-Lie cases.

An important example of an augmented Lie rack is a Lie group crossed module. First, we recall the definition of Lie group crossed modules.

\begin{de.}

A Lie group crossed module is a 4-tuple $(M,N,\mu,\eta)$ where $M,N$ is a Lie group, $\mu:M\rightarrow N$ is a morphism of Lie groups and $\eta:N\rightarrow \mathrm{Aut}(M)$ is a homomorphism defining a smooth action $\hat{\eta}:N\times M\rightarrow M$. They satisfy the following conditions for all $n\in N,m,m'\in M$.

\begin{enumerate}

\item$\mu(\eta(n)(m))=n\mu(m)n^{-1}$,

\item$\eta(\mu(m))(m')=mm'm^{-1}$.

\end{enumerate}

\end{de.}

With an $N$-action $\hat{\eta}$, $(N,M,\mu)$ becomes an augmented Lie rack. More generally, a group crossed module is an instance of augmented racks.

We can relax the condition in a similar way that augmented Lie racks are relaxed.

\begin{de.}

A relaxed Lie group crossed module is a 5-tuple $(M,N,N',\mu,\eta)$ where $M,N$ is a Lie group, $\mu:M\rightarrow N$ is a morphism of Lie groups,$N'\subset N$ is a Lie subgroup of $N$ including $\mu(M)$ and $\eta:N\rightarrow \mathrm{Aut}(M)$ defines a smooth action $\hat{\eta}:N\times M\rightarrow M$. They satisfy the two conditions of a Lie group crossed module for all $n\in N',m,m'\in M$.

\end{de.}

When $N'=N$, it reduces to an ordinary Lie group crossed module. With an $N$-action $\hat{\eta}$, $(N,M,\mu,N')$ becomes a relaxed augmented Lie rack and hence $(N,M,\mu)$ becomes a Lie group-rack triple. Note that we can define a relaxed (non-Lie) group crossed module and treat it as an example of relaxed augmented racks.

It is known that, for a Lie group crossed module $(M,N,\mu,\eta)$, $(T_{e_{M}}M,T_{e_{N}}N,T_{e_{M}}\mu,T_{e_{N}}\eta)$ becomes a Lie algebra crossed module. This fact also holds for the relaxed case.

We show that the general tangent structure of a relaxed augmented Lie rack is a relaxed augmented Leibniz algebra.

\begin{pr.}
Let $(G,(X,e),\Theta:X\rightarrow G,G')$ be a relaxed augmented Lie rack, $\mathfrak{g},\mathfrak{g'}$ be the Lie algebra of $G,G'$, and $V:=T_{e}X$. We give $V$ a structure of $\mathfrak{g}$-module via the differential at $e\in X$ of the Lie group action of $G$ on $X$. Explicitly, we define
\begin{equation}
a\cdot v=\left.\frac{d}{d t}T_{e}L_{a(t)}(v)\right|_{t=0}
\end{equation}
for all $a\in\mathfrak{g}$ and $v\in V$, where $t\mapsto a(t)\in G$ is any smooth curve on $t\in[-1,1]$ satisfying $a(0)=e_{G}$ and $da/dt(0)=a\in\mathfrak{g}$ and $L_{a(t)}$ is the Lie group action of $G$ on $X$. Let $\theta=d\Theta:V\rightarrow\mathfrak{g}$ be a differential of $\Theta$ at $e\in X$. Then $(\mathfrak{g},V,\theta,\mathfrak{g'})$ is a relaxed augmented Leibniz algebra. We call this 4-tuple the tangent relaxed augmented Leibniz algebra of $(G,(X,e),\Theta:X\rightarrow G,G')$.

Moreover, let $(G_{i},X_{i},\Theta_{i})(i=1,2)$ be the underlying Lie group-rack triples of a relaxed augmented Lie rack $(G_{i},(X_{i},e_{i}),\Theta_{i}:X_{i}\rightarrow G_{i},G'_{i})$ and $(\Phi,\Psi):(G_{1},X_{1},\Theta_{1})\rightarrow(G_{2},X_{2},\Theta_{2})$ a morphism of Lie group-rack triples. Then, $(\varphi=T_{e_{G_{1}}}\Phi,\psi=T_{e_{X_{1}}}\Psi):(\mathfrak{g}_{1},V_{1},\theta_{1})\rightarrow(\mathfrak{g}_{2},V_{2},\theta_{2})$ is a morphism of Lie-Leibniz triples where $(\mathfrak{g}_{i},V_{i},\theta_{i},\mathfrak{g}'_{i})$ are the tangent relaxed augmented Leibniz algebra of $(G_{i},(X_{i},e_{i}),\Theta_{i}:X_{i}\rightarrow G_{i},G'_{i})$.

\end{pr.}

\begin{proof}

Let $t\mapsto v(t)\in X$ be a smooth curve on $t\in[-1,1]$ satisfying $v(0)=e$ and $dv/dt(0)=v\in V$. Choose $a(t)$ as $a(t)\in G'$ when $t\in (-\epsilon,\epsilon)$ for small $\epsilon$. In this case, $a=da/dt(0)\in\mathfrak{g}'$. Then we compute

\begin{align*}
\theta(a\cdot v)&=d\Theta\left(\left.\frac{\partial}{\partial t}T_{e}L_{a(t)}(v)\right|_{t=0}\right)\\
&=\left.\frac{\partial}{\partial t_{2}}\frac{\partial}{\partial t_{1}}\Theta\left(L_{a(t_{1})}(v(t_{2}))\right)\right|_{t_{1},t_{2}=0}\\
&=\left.\frac{\partial}{\partial t_{2}}\frac{\partial}{\partial t_{1}}a(t_{1})\Theta(v(t_{2}))a(t_{1})^{-1}\right|_{t_{1},t_{2}=0}\\
&=\left[\left.\frac{\partial}{\partial t_{1}}a(t_{1})\right|_{t_{1}=0},\left.\frac{\partial}{\partial t_{2}}\Theta(v(t_{2}))\right|_{t_{2}=0}\right]_{\mathfrak{g}}=[a,\theta(v)]_{\mathfrak{g}}.
\end{align*}
The second equality follows from $L_{a(t_{1})}(v(0))=e$ and $\left.\frac{\partial}{\partial t_{2}}L_{a(t_{1})}(v(t_{2}))\right|_{t_{2}=0}=T_{e}L_{a(t_{1})}(v)$.

For the second statement, let $t\mapsto v_{1}(t)\in X_{1}$ be a smooth curve on $t\in[-1,1]$ satisfying $v(0)=e_{X_{1}}$ and $dv/dt(0)=v_{1}\in V_{1}$. Then $\Psi(v_{1}(0))=e_{X_{2}}$ and $\frac{d}{dt}\Psi(v_{1}(t))|_{t=0}=\psi(v_{1})$. We compute
\begin{align*}
\theta_{2}(\psi(v_{1}))&=T_{e_{X_{2}}}\Theta_{2}(T_{e_{X_{1}}}\Psi(v_{1}))\\
&=\left.\frac{d}{dt}\Theta_{2}(\Psi(v_{1}(t)))\right|_{t=0}\\
&=\left.\frac{d}{dt}\Phi(\Theta_{1}(v_{1}(t)))\right|_{t=0}\\
&=T_{e_{X_{2}}}\Phi(T_{e_{X_{1}}}\Theta_{1}(v_{1}))=\varphi(\theta_{1}(v_{1})),
\end{align*}

\begin{align*}
\psi(a\cdot v_{1})&=T_{e_{X_{1}}}\Psi\left(\left.\frac{d}{dt}T_{e}L_{a(t)}(v_{1})\right|_{t=0}\right)\\
&=\left.\frac{\partial}{\partial t_{1}}\frac{\partial}{\partial t_{2}}\Psi(L_{a(t_{1})}v_{1}(t_{2}))\right|_{t_{1},t_{2}=0}\\
&=\left.\frac{\partial}{\partial t_{1}}\frac{\partial}{\partial t_{2}}\Phi(a(t_{1}))\cdot\Psi(v_{1}(t_{2}))\right|_{t_{1},t_{2}=0}\\
&=\left.\frac{d}{dt}T_{e}L_{\Phi(a(t))}(\psi(v_{1}))\right|_{t=0}=\varphi(a)\cdot\psi(v_{1}).
\end{align*}

\end{proof}

Define a bracket $[,]$ on $V$ by
\begin{align*}
[v,u]&=\theta(v)\cdot u\\
&=\left.\frac{\partial}{\partial t}T_{e}L_{\Theta(v(t))}(u)\right|_{t=0}\\
&=\left.\frac{\partial}{\partial t}T_{e}L_{v(t)}(u)\right|_{t=0}
\end{align*}
where $L_{v(t)}:=v(t)\triangleright-:X\rightarrow X$. Then $(V,[,])$ becomes a Leibniz algebra. The last equality is the original definition of the tangent Leibniz algebra of a given Lie rack\cite{kinyon2007leibniz}.




For a general curve $a(t)$ in $G$ such that $a(0)=e_{G}$ and $da/dt(0)=a$, we have
\begin{align*}
&\ \left.\frac{\partial}{\partial t_{2}}\frac{\partial}{\partial t_{1}}(a(t_{1})\Theta)(v(t_{2}))\right|_{t_{1},t_{2}=0}\\
&=\left.\frac{\partial}{\partial t_{2}}\frac{\partial}{\partial t_{1}}a(t_{1})\Theta(v(t_{2}))a(t_{1})^{-1}\left(\Theta\left(L_{a(t_{1})}(v(t_{2}))\right)\right)^{-1}\right|_{t_{1},t_{2}=0}\\
&=\left[\left.\frac{\partial}{\partial t_{1}}a(t_{1})\right|_{t_{1}=0},\left.\frac{\partial}{\partial t_{2}}\Theta(v(t_{2}))\right|_{t_{2}=0}\right]_{\mathfrak{g}}-d\Theta\left(\left.\frac{\partial}{\partial t}T_{e}L_{a(t)}(v)\right|_{t=0}\right)\\
&=[a,\theta(v)]_{\mathfrak{g}}-\theta(av)=(a\theta)(v).
\end{align*}

Next, we define the local version of relaxed augmented Lie racks. The local version of (augmented) Lie racks was introduced in \cite{covez2010integration}.

\begin{de.}
A local rack is a set $X$ endowed with a binary operation $\triangleright:\Omega\rightarrow X$ where $\Omega\subset X\times X$ which satisfies the following conditions.

1.If $(x,y),(y,z),(x,z),(x,y\triangleright z),(x\triangleright y,x\triangleright z)\in\Omega$, then self-distributivity $x\triangleright(y\triangleright z)=(x\triangleright y)\triangleright (x\triangleright z)$ holds.

2.If $(x,y),(x,z)\in\Omega$ and $x\triangleright y=x\triangleright z$, then $y=z$.

We say that a local rack is pointed if there exists $e\in X$ such that $e\triangleright x=x$ and $x\triangleright e=e$ for all $x\in X$. A local Lie rack is a pointed local rack $(X,\Omega,e)$ where $X$ is a smooth manifold, $\Omega$ is an open subset of $X\times X$ and $\triangleright:\Omega\rightarrow X$ is smooth.
\end{de.}

The local version of a Lie group action on a manifold is defined as follows. 

\begin{de.}
Let $G$ be a group with a unit element $e_{G}$. A local $G$-set is a set $X$ with a map $q:\Omega\rightarrow X$ where $\Omega\subset G\times X$ which satisfies the following conditions.

\begin{enumerate}

\item If $(h,x),(gh,x),(g,q(h,x))\in\Omega$, then $q(g,q(h,x))=q(gh,x)$.

\item $(e_{G},x)\in\Omega$ and $q(e_{G},x)=x$ for all $x\in X$.

\end{enumerate}

Assume that $G$ is a Lie group. A local smooth $G$-set is a smooth manifold $X$ with a local $G$-set structure where $\Omega$ is an open subset of $G\times X$ and $q:\Omega\rightarrow X$ is smooth. A fixed point of a local smooth $G$-set is $e_{X}\in X$ such that $(g,e_{X})\in\Omega$ and $q(g,e_{X})=e_{X}$ for all $g\in G$.
\end{de.}

With the above definition, we can define a local relaxed augmented Lie rack.

\begin{de.}

Let $G$ be a Lie group, and $X$ be a local smooth $G$-set with a map $q:\Omega\rightarrow X$, which has a fixed point $e_{X}\in X$. We say that $X$ together with a smooth map $\Theta:X\rightarrow G$ and a Lie subgroup $G'\subset G$ including $\mathrm{Im}(\Theta)$ is a local relaxed augmented Lie rack when it satisfies
\begin{equation}
\Theta(e_{X})=e_{G}
\end{equation}
and
\begin{equation}
\Theta(q(g,x))=g\Theta(x)g^{-1}
\end{equation}
for all $(g,x)\in G'\times X\cap\Omega$. 
    
\end{de.}

For a local relaxed augmented Lie rack $(G,(X,e_{X}),q:\Omega\rightarrow X,,\Theta,G')$, $(X,\Omega',e_{X})$ becomes a local Lie rack with
\begin{equation}
x\triangleright y:=q(\Theta(x),y)
\end{equation}
for all $(x,y)\in\Omega'$ where
\begin{equation}
\Omega':=\{(x,y)\in X\times X|(\Theta(x),y)\in\Omega\}.
\end{equation}

Let $\mathfrak{g}$ and $\mathfrak{g}'$ be the tangent Lie algebra of $G$ and $G'$, $V:=T_{e_{X}}X$ be the tangent space of $X$, and $dq:\mathfrak{g}\times V\rightarrow V$ and $\theta=d\Theta:V\rightarrow\mathfrak{g}$ be the differential of $q$ and $p$. Then $(\mathfrak{g},V,\theta,\mathfrak{g}')$ becomes a relaxed augmented Leibniz algebra and $(\mathfrak{g},V,\theta)$ becomes a Lie-Leibniz triple. We call this triple the tangent Lie-Leibniz triple of $(G,(X,e_{X}),q:\Omega\rightarrow X,\Theta,G')$.

\section{local integration procedure}

In this section, we give the local integration procedure of real finite-dimensional Lie-Leibniz triples. We integrate a relaxed augmented Leibniz algebra into a local relaxed augmented Lie rack by generalizing the integration for the strict case in \cite{bordemann2017global}. 

Let $(\mathfrak{g},V,\theta,\mathfrak{h})$ be the relaxed augmented Leibniz algebra associated with a Lie-Leibniz triple $(\mathfrak{g},V,\theta)$ where both $\mathfrak{g}$ and $V$ is finite-dimensional. We set $\mathfrak{g}':=\mathrm{Im}(\theta),\mathfrak{z}:=\mathrm{Ker}(\theta)$. Let $G,G',H$ be the connected, simply connected Lie groups associated to $\mathfrak{g},\mathfrak{g}',\mathfrak{h}$. Let $\tilde{G}',\tilde{H}\subset G$ be the connected Lie subgroup of $G$ whose tangent Lie algebra is $\mathfrak{g}',\mathfrak{h}$.

First, we make some preparations to generalize the integration method for the strict case.

Consider the adjoint representation of $H$ on $\mathfrak{h}$. Since this representation preserves ideals and $\mathfrak{g}'$ is an ideal of the Lie algebra $\mathfrak{h}$, for all $h\in H$, $Ad_{h}:\mathfrak{h}\rightarrow\mathfrak{h}$ defines an automorphism of $\mathfrak{g}'$. The correspondence $H\ni h\mapsto Ad_{h}\in \mathrm{Aut}(\mathfrak{g}')$ defines a Lie group homomorphism between $H$ and $\mathrm{Aut}_{0}(\mathfrak{g}')$, the connected component of the Lie group $\mathrm{Aut}(\mathfrak{g}')$ including the identity. As a Lie group, $\mathrm{Aut}_{0}(\mathfrak{g}')$ is isomorphic to $\mathrm{Aut}_{0}(G')$. Therefore, there is a unique Lie group homomorphism
\begin{equation}
I':H\rightarrow\mathrm{Aut}_{0}(G'), h\mapsto(g'\mapsto I'_{h}(g'))
\end{equation}
such that $T_{e_{G}}I'_{h}=Ad_{h}$ for all $h\in H$. 

We use the following lemma.

\begin{le.}[{\cite[Lemma 3.1]{bordemann2017global}}]

There is a smooth map $s:G'\rightarrow\mathfrak{g}'$ that has the following properties:
\begin{equation}
s(e_{G'})=0,
\end{equation}
\begin{equation}
T_{e_{G'}}s=id_{\mathfrak{g'}},
\end{equation}
\begin{equation}
\forall h\in H,\forall g'\in G':s(I'_{h}(g'))=\mathrm{Ad}_{h}(s(g')).
\end{equation}

\end{le.}

(For the proof and construction of $s$, see the proof of Lemma 3.1 of \cite{bordemann2017global}.)

Next, we consider the immersion of $G'$ into $H$ and of $H$ into $G$. Let $\iota_{G'}:G'\hookrightarrow H$ and $\iota_{H}:H\hookrightarrow G$ be the immersions of Lie groups derived from the inclusions $\mathfrak{g}'\hookrightarrow\mathfrak{h}$ and $\mathfrak{h}\hookrightarrow\mathfrak{g}$ of Lie algebras. We have
\begin{equation}
(\iota_{H}\circ\iota_{G'}\circ I'_{h}(g'))=\iota_{H}(h)(\iota_{H}\circ\iota_{G'}(g'))\iota_{H}(h)^{-1}
\end{equation}
for all $h\in H$ and $g'\in G'$.

Let $\rho:G\times V\rightarrow V$ be the unique representation of $G$ on $V$ such that for all $\xi\in\mathfrak{g}$ and $v\in V$
\begin{equation}
\left.\frac{d}{dt}\rho_{exp(t\xi)}(v)\right|_{t=0}=\xi\cdot v.
\end{equation}

The map
\begin{equation}
H\times V\ni(h,v)\mapsto\rho_{\iota_{H}(h)}(v)\in V
\end{equation}
gives the representation of $H$ on $V$ whose differential representation is the restricted representation of the $\mathfrak{g}$-module structure on $V$ to $\mathfrak{h}$.

Since $\theta(\xi v)=ad_{\xi}(\theta(v))$ when $\xi\in\mathfrak{h}$ and $v\in V$, we have
\begin{equation}
\theta(\rho_{\iota_{H}(h)}(v))=Ad_{h}(\theta(v))
\end{equation}
 for all $h\in H$ and $v\in V$.

The following theorem states that, for any real finite-dimensional Lie-Leibniz triple, we can construct a local relaxed augmented Lie rack which integrates the triple. 

\begin{th.}

Consider $\theta:V\rightarrow\mathfrak{g}'$ as an affine bundle over $\mathfrak{g'}$ with typical fiber $\mathfrak{z}$.

Define the pull-back fiber bundle
\begin{equation}
M:=s^{*}V=\{(v,g')\in V\times G'|\theta(v)=s(g')\}\rightarrow G'
\end{equation}
over $G'$ and let $(0,e_{G})\in M$ be the distinguished point. Let $\phi:M\rightarrow G'$ be the map $M\ni(v,g')\mapsto g'\in G'$. Then there is an open submanifold $M_{U}$ of $M$, an open subset $\Omega\subset G\times M_{U}$ and a local smooth $G$-action $q:\Omega\rightarrow M_{U}$ such that $(G,(M_{U},(0,e_{G'})),q,\Phi:M_{U}\rightarrow G,\tilde{H})$, where we write the map $\iota_{H}\circ\iota_{G'}\circ(\phi|_{M_{U}}):M_{U}\rightarrow G$ as $\Phi$, is a local relaxed augmented Lie rack such that the tangent structure becomes the relaxed augmented Leibniz algebra $(\mathfrak{g},V,\theta,\mathfrak{h})$. When the associated Lie-Leibniz triple is strict, the local $G$-action can be extended to the Lie group action on $M$ and $(G,(M,(0,e_{G'})),\Phi)$ becomes an augmented Lie rack.

\end{th.}

\begin{proof}

First, consider the Lie algebra action of $\mathfrak{g}$ on the tangent space of $M$ at $(0,e_{G'})$. Due to the fact that $T_{e_{G'}}s=id_{\mathfrak{g}}$, we have

\begin{equation}
T_{(0,e_{G'})}M=\{(v,\zeta)\in T_{0}V\times T_{e_{G'}}G'|\theta(v)=\zeta\}.
\end{equation}

We define the $\mathfrak{g}$-module structure on $T_{(0,e_{G'})}M$ via
\begin{equation}
a(v,\zeta):=(av,\theta(av))
\end{equation}
for all $a\in\mathfrak{g}$. When $a\in\mathfrak{h}$, we get $a(v,\zeta)=(av,[a,\theta(v)]_{\mathfrak{g}})$.

Next, we consider the local Lie group action of $G$ on the open neighborhood of $(0,e_{G'})$ in $M$.

Due to the fact $T_{e_{G'}}s=id_{\mathfrak{g}}$ and the inverse function theorem, there is an open neighborhood $U\subset G'$ of $e_{G'}$ such that $s|_{U}:U\rightarrow s(U)$ is a diffeomorphism.
Let $s^{-1}:s(U)\rightarrow U$ be the inverse diffeomorphism of $s$.

We define
\begin{equation}
M_{U}:=\{(v,u)\in V\times U | \theta(v) = s (u)\}\subset M.
\end{equation}
Note that $u$ is uniquely determined by $u=s^{-1}(\theta(v))$ on $M_{U}$.

First, we show that $M_{U}$ is a local smooth $G$-set. We define an open subset of $G\times M_{U}$ as
\begin{equation}
\Omega:=\{(g,v,u)\in G\times M_{U}|\theta(\rho_{g}(v))\in s(U)\}.
\end{equation}

We define a map defined on $\Omega$ with values in $M_{U}$ by
\begin{equation}
\Omega\ni (g,v,u)\mapsto g(v,u)=(\rho_{g}(v),s^{-1}(\theta(\rho_{g}(v)))).
\end{equation}

We can check that $M_{U}$ with a map $q:\Omega\rightarrow M_{U}$ satisfies the conditions of a local smooth $G$-set as follows. 

1.If $(g_{1},v),(g_{1}g_{2},v),(g_{1},q(g_{2},v))\in\Omega$, then 
\begin{align*}
q(g_{1},q(g_{2},v,u))&=q(g_{1},\rho_{g_{2}}(v),s^{-1}(\theta(\rho_{g_{2}}(v))))\\
&=(\rho_{g_{1}}(\rho_{g_{2}}(v)),s^{-1}(\theta(\rho_{g_{1}}(\rho_{g_{2}}(v)))))\\
&=(\rho_{g_{1}g_{2}}(v),s^{-1}(\theta(\rho_{g_{1}g_{2}}(v))))\\
&=q(g_{1}g_{2},v,u).  
\end{align*}

2.$(e_{G},v,u)\in\Omega$ and 
\begin{align*}
q(e_{G},v,u)&=(\rho_{e_{G}}(v),s^{-1}(\theta(\rho_{e_{G}}(v))))\\
&=(v,s^{-1}(\theta(v)))=(v,u)\\
\end{align*}
for all $v\in M_{U}$. 

$(0,e_{G'})\in M_{U}$ is a fixed point of the local smooth $G$-action. Indeed, $\theta(\rho_{g}(0))=0\in s(U)$ and
\begin{align*}
q(g,0,e_{G'})&=(\rho_{g}(0),s^{-1}(\theta(\rho_{g}(0))))\\
&=(0,e_{G'})\\
\end{align*}
for all $v\in M_{U}$.

Consider the local smooth action of the Lie subgroup $\tilde{H}$ of $G$. Note that there is a smooth action of $H$ on $M$ via
\begin{equation}\label{h}
M\ni(v,g')\mapsto (\rho_{\iota_{H}(h)}(v),I'_{h}(g'))\in M
\end{equation}
for all $h\in H$. This action is well-defined because 
\begin{equation}
\theta(\rho_{\iota_{H}(h)}(v))=Ad_{h}(\theta(v))=Ad_{h}(s(g'))=s(I'_{h}(g')).
\end{equation}

When $\tilde{h}\in \tilde{H}$, the local smooth $G$-set structure on $M_{U}$ coincides with the action of $H$. That is, there is an element $h\in H$ such that $\iota_{H}(h)=\tilde{h}$ and
\begin{align*}
\tilde{h}(v,u)&=(\rho_{\tilde{h}}(v),s^{-1}(\theta(\rho_{\tilde{h}}(v))))\\&=(\rho_{\tilde{h}}(v),s^{-1}(Ad_{h}\theta(v)))\\
&=(\rho_{\tilde{h}}(v),I'_{h}(u)).
\end{align*}
We can check this claim as follows. Since $\iota_{H}(H)=\tilde{H}$, there is an element $h\in H$ such that $\iota_{H}(h)=\tilde{h}$ and $s(I'_{h}(u))=Ad_{h}\theta(v)=\theta(\rho_{\tilde{h}}(v))$. When $(\tilde{h},v,u)\in\Omega\cap \tilde{H}\times M_{U}$, $\theta(\rho_{\tilde{h}}(v))\in s(U)$ and $s^{-1}(\theta(\rho_{\tilde{h}}(v)))=I'_{h}(u)$ is uniquely determined.

Next, we consider the structure of a local relaxed augmented Lie rack on $M_{U}$. 

For any $(v,u)\in M_{U}$ and $h\in \tilde{H}$, there is an element $h\in H$ such that $\iota_{H}(h)=\tilde{h}$ and
\begin{align*}
\Phi(q(\tilde{h},v,u))&=\Phi(\rho_{\tilde{h}}(v),s^{-1}(\theta(\rho_{\tilde{h}}(v))))\\
&=\Phi(\rho_{\tilde{h}}(v),s^{-1}(Ad_{h}(\theta(v))))\\
&=\Phi(\rho_{\tilde{h}}(v),I'_{h}(s^{-1}(\theta(v))))\\
&=\iota_{H}\circ\iota_{G'}\circ I'_{h}(u)\\
&=\iota_{H}(h)(\iota_{H}\circ\iota_{G'}(u))\iota_{H}(h^{-1})\\
&=\tilde{h}(\Phi(v,u))\tilde{h}^{-1}.
\end{align*}
This shows that $((M_{U},(0,e_{G'})),G,q,\Phi,\tilde{H})$ is a local relaxed augmented Lie rack.

We consider the tangent structure of the local relaxed augmented Lie rack. By the definition,

\begin{equation}
T_{(0,e_{G'})}M_{U}=T_{(0,e_{G'})}M=\{(v,\zeta)\in T_{0}V\times T_{e_{G'}}G'|\theta(v)=\zeta\}.
\end{equation}

Let $dq:\mathfrak{g}\times T_{(0,e_{G'})}M_{U}\rightarrow T_{(0,e_{G'})}M_{U}$ and $d\Phi:T_{(0,e_{G'})}M_{U}\rightarrow\mathfrak{g}$ be the differential of $q$ and $d\Phi$. Then $(\mathfrak{g},T_{(0,e_{G'})}M_{U},d\Phi,\mathfrak{h})$ becomes a relaxed augmented Leibniz algebra. We construct an isomorphism between two Lie-Leibniz triples $(\mathfrak{g},V,\theta)$ and $(\mathfrak{g},T_{(0,e_{G'})}M_{U},d\Phi)$ as a pair $(id_{\mathfrak{g}},\chi)$ where
\begin{equation}
\chi:T_{(0,e_{G'})}M_{U}\ni(v,\theta(v))\mapsto v\in V.
\end{equation}

We can check that $(id_{\mathfrak{g}},\chi)$ is an isomorphism as follows. We choose a smooth curve $[-1,1]\ni t\mapsto u(t)\in G'$ such that $u(0)=e_{G'}$ and $\frac{d}{dt}u(t)|_{t=0}=\theta(v)$. Then we have

\begin{align*}
d\Phi(v,\theta(v))&=\left.\frac{\partial}{\partial t_{2}}\frac{\partial}{\partial t_{1}}\Phi(t_{1}v,u(t_{2})) \right|_{t_{1},t_{2}=0}\\
   &=\left.\frac{d}{dt}u(t) \right|_{t=0}=\theta(v)=\theta(\chi(v,\theta(v))),
\end{align*}

\begin{align*}
\chi(dq(a,v,\theta(v)))&=\chi\left(\left.\frac{\partial}{\partial t_{2}}\frac{\partial}{\partial t_{1}}(\rho_{a(t_{2})}(t_{1}v),s^{-1}(\theta(\rho_{a(t_{2})}(t_{1}v))))\right|_{t_{1},t_{2}=0}\right)\\
&=\chi\left(\left.\frac{d}{d t}(a\cdot t_{1}v,s^{-1}(\theta(a\cdot t_{1}v)))\right|_{t=0}\right)\\
&=\chi(a\cdot v,\theta(a\cdot v))=a\cdot v=a\cdot\chi(v,\theta(v)).
\end{align*}

These calculations show that $(id_{\mathfrak{g}},\chi)$ is a morphism of Lie-Leibniz triples. It is clear that both $id_{\mathfrak{g}}$ and $\chi$ are isomorphisms. Thus, $(id_{\mathfrak{g}},\chi)$ is an isomorphism of Lie-Leibniz triples.

When the tangent Lie-Leibniz triple is strict, then we can take $H=\tilde{H}=G$ and $\iota_{H}=id_{G}$, and we can extend the local $G=H$ action to the action on the whole manifold $M$ via 
\begin{equation}
M\ni(v,g')\mapsto (\rho_{h}(v),I'_{h}(g'))\in M
\end{equation}
for all $h\in H$ and $((M,(0,e_{G'})),G=H,q,\Phi)$ becomes an augmented Lie rack.

\end{proof}

Note that $(H,M,\iota_{G'}\circ\phi:M\rightarrow H)$ with the action of $H$ on $M$ (\ref{h}) has the structure of an augmented Lie rack integrating the augmented Leibniz algebra $(\mathfrak{h},V,\theta)$, where $V$ is seen as an $\mathfrak{h}$-module and $\theta$ as a map $\theta:V\rightarrow\mathfrak{h}$. What we have shown is that $M$ also captures the integrated structure of the Lie-Leibniz triple $(\mathfrak{g},V,\theta)$, where $V$ is seen as an $\mathfrak{g}$-module and $\theta$ as a map $\theta:V\rightarrow\mathfrak{g}$, locally.

\section{future directions}

As an integration of Lie-Leibniz triples, the integration process introduced in this paper is local. The existence of the global integration process of general Lie-Leibniz triples is unknown.

The Lie algebras of many important Lie-Leibniz triples are infinite-dimensional. An infinite-dimensional Lie group-rack triple and the integration procedure for an infinite-dimensional Lie-Leibniz triple is yet to be understood.

An important problem is whether we can obtain the global structure associated to the tensor hierarchy algebra\cite{Palmkvist:2013vya} from a given Lie-Leibniz triple. Like a Lie-Leibniz triple gives rise to a dgLa, a Lie group-rack triple would give rise to a global object integrating the dgLa, such as a differential-graded Lie group\cite{jubin2022differential}.

A Lie group-rack triple is a global counterpart of the Lie-Leibniz triple. Thus, it would be a useful tool for analyzing global aspects of tensor hierarchies. In particular, we would be able to develop the theory of "principal bundles of Lie racks" to consider the analog of ordinary "principal bundles of Lie groups" for Leibniz gauge theories\cite{Bonezzi:2019ygf}.

$\,$

$\,$

\bibliography{myref}
\bibliographystyle{unsrt}

\end{document}